\pdfoutput=1
\newif\ifSpringer

\Springerfalse 

\ifSpringer

\RequirePackage{fix-cm}
\documentclass[smallextended]{svjour3}
\usepackage{etoolbox}
\AtEndEnvironment{proof}{\phantom{}\qed}

\usepackage{amsmath,amssymb,graphicx}

\else

\documentclass[11pt]{scrartcl}
\usepackage{amsmath,amssymb,amsthm,graphicx}
\newtheorem{theorem}{Theorem}[section]

\newtheorem{example}[theorem]{Example}
\fi

\usepackage[utf8x]{inputenc}
\usepackage[T1]{fontenc}
\usepackage{lmodern}
\usepackage[english]{babel}

\usepackage[draft=false,kerning=true]{microtype}

\usepackage{todonotes}
\usepackage{color}
\usepackage{colortbl}
\usepackage{mathtools}
\usepackage{hyperref}
\usepackage{longtable}
\usepackage{geometry}

\usepackage[font=small]{subcaption}

\usepackage{scalerel}
\usepackage{tikz}
\usetikzlibrary{svg.path}

\usepackage{pgfplots}
\pgfplotsset{compat=1.12}

\usepackage[ruled,vlined,linesnumbered]{algorithm2e}

\usepackage[normalem]{ulem}

\DeclareMathOperator{\diag}{diag}
\DeclarePairedDelimiter{\abs}{\lvert}{\rvert}
\DeclarePairedDelimiter\set\{\}
\def\CPP{{C\nolinebreak[4]\hspace{-.05em}\raisebox{.4ex}{\tiny ++}}}
\newcommand{\R}{\mathbb{R}}

\newcommand{\elli}[1]{#1}
\newcommand{\melanie}[1]{#1}

\usepackage{bbm} 
\DeclareMathOperator{\STAB}{STAB}

\DeclareMathOperator{\ThetaBody}{TH}
\DeclareMathOperator{\conv}{conv}

\DeclareMathOperator{\st}{s.t.}

\newcommand{\Sym}[1]{{\mathcal S}_{#1}}
\newcommand{\simplex}{{\Delta}}
\newcommand{\allones}[1]{\mathbbm{1}_{#1}}
\newcommand{\transposedVec}{T}
\newcommand{\norm}[1]{\left\lVert #1 \right\rVert}

\newcommand{\GI}{G_{I}}
\newcommand{\kI}{k_{I}}
\newcommand{\bI}{b_{I}}

\newcommand{\tI}{t_{I}}

\newcommand{\ofvSDPescC}{z_J^{C}(G)}
\newcommand{\ofvSDPescF}{z_J^{F}(G)}
\newcommand{\ofvSDPescH}{z_J^{H}(G)}
\newcommand{\XI}{X_{I}}
\newcommand{\SIi}{S^{I}_{i}}
\newcommand{\sIi}[1]{s^{I}_{#1}}
\newcommand{\SIiPrime}{S^{I}_{i'}}
\newcommand{\zESCJ}{z_{J}(G)}

\newcommand{\JbNew}{J_{\nrESCNew}}

\newcommand{\XOptC}{X^{C\ast}}

\newcommand{\GZeroi}[1]{G^{0}_{#1}}
\newcommand{\VZeroi}[1]{V^{0}_{#1}}
\newcommand{\EZero}{E^{0}}

\newcommand{\Fik}[1]{F^{#1}_{i}}
\newcommand{\fik}[1]{f^{#1}_{i}}
\newcommand{\nrFacetsk}[1]{r_{#1}}
\newcommand{\epsF}{\varepsilon_{F}}
\newcommand{\XIOpt}{X_{I}^{\ast}}
\newcommand{\indicesVFIPrime}{\mathcal{V}_{I}'}
\newcommand{\indicesVFI}{\mathcal{V}_{I}}
\newcommand{\FiPrimek}[1]{F^{#1}_{i'}}
\newcommand{\fiPrimek}[1]{f^{#1}_{i'}}
\newcommand{\XOptF}{X^{F\ast}}
\newcommand{\XOptFI}{\XOptF_{I}}

\newcommand{\projXI}{P_{I}}
\newcommand{\QI}{Q_{I}}
\newcommand{\HI}{H_{I}}
\newcommand{\hI}{h_{I}}
\newcommand{\XOptH}{X^{H\ast}}
\newcommand{\XOptHI}{\XOptH_{I}}
\newcommand{\XHyp}{X^{\ast}}

\newcommand{\Xtilde}{\tilde{X}}
\newcommand{\XtildeI}{\Xtilde_{I}}

\newcommand{\XOpt}{X^{\ast}}
\newcommand{\XOptI}{X^{\ast}_I}
\newcommand{\lambdaI}{\lambda_{I}}
\newcommand{\lambdaIi}[1]{[\lambdaI]_{#1}}
\newcommand{\nrESC}{|J|}
\newcommand{\nrESCNew}{q}

\newcommand{\mytag}[1]{(\hypertarget{#1}{#1})}
\newcommand{\myref}[1]{(\hyperlink{#1}{#1})}

\newcommand{\nF}{\cellcolor{red!25}}
\newcommand{\nO}{\cellcolor{yellow!45}}

\definecolor{orcidlogocol}{HTML}{A6CE39}
\tikzset{
	orcidlogo/.pic={
		\fill[orcidlogocol] svg{M256,128c0,70.7-57.3,128-128,128C57.3,256,0,198.7,0,128C0,57.3,57.3,0,128,0C198.7,0,256,57.3,256,128z};
		\fill[white] svg{M86.3,186.2H70.9V79.1h15.4v48.4V186.2z}
		svg{M108.9,79.1h41.6c39.6,0,57,28.3,57,53.6c0,27.5-21.5,53.6-56.8,53.6h-41.8V79.1z M124.3,172.4h24.5c34.9,0,42.9-26.5,42.9-39.7c0-21.5-13.7-39.7-43.7-39.7h-23.7V172.4z}
		svg{M88.7,56.8c0,5.5-4.5,10.1-10.1,10.1c-5.6,0-10.1-4.6-10.1-10.1c0-5.6,4.5-10.1,10.1-10.1C84.2,46.7,88.7,51.3,88.7,56.8z};
	}
}

\newcommand\orcidicon[1]{\href{https://orcid.org/#1}{\mbox{\scalerel*{
				\begin{tikzpicture}[yscale=-1,transform shape]
				\pic{orcidlogo};
				\end{tikzpicture}
			}{|}}}}

\providecommand{\keywords}[1]
{
	\small	
	\textbf{\textit{Keywords---}} #1
}

\title{An SDP-Based Approach for Computing\\
 	the Stability Number
 	of a Graph%
        \ifSpringer
        \else
        \footnote{This project has received funding from the Austrian 
          Science Fund (FWF): I\,3199-N31.
          Moreover, the second author has received funding from the Austrian 
          Science Fund (FWF): DOC~78.
      \elli{Furthermore, we thank two anonymous referees for their valuable 
      input.}}%
        \fi
}

\ifSpringer

\author{Elisabeth Gaar \and Melanie Siebenhofer \and Angelika Wiegele}
\authorrunning{E. Gaar, M. Siebenhofer, A. Wiegele} 
\institute{
        E. Gaar \at
        Institute of Production and Logistics Management, 
        Johannes Kepler University Linz, 
        Altenberger Stra{\ss} 69, 4040 Linz, Austria,
        ORCiD: 0000-0002-1643-6066, \email{elisabeth.gaar@jku.at}
  	\and
        M. Siebenhofer \at
	Institut f\"ur Mathematik, Alpen-Adria-Universit\"at Klagenfurt, Universit\"atsstra{\ss}e 65-67, 9020 Klagenfurt, Austria,
	ORCiD: 0000-0002-9101-834X, \email{melanie.siebenhofer@aau.at}
        \and
        A. Wiegele \at
        Institut f\"ur Mathematik, Alpen-Adria-Universit\"at Klagenfurt, Universit\"atsstra{\ss}e 65-67, 9020 Klagenfurt, Austria,
	ORCiD: 0000-0003-1670-7951, \email{angelika.wiegele@aau.at}  
}

\else

\author{Elisabeth Gaar \orcidicon{0000-0002-1643-6066}, Melanie Siebenhofer \orcidicon{0000-0002-9101-834X} and Angelika Wiegele \orcidicon{0000-0003-1670-7951}}

\date{\today}

\fi

\begin{document}

\maketitle

\begin{abstract}
Finding the stability number of a graph, i.e., the maximum number of vertices 
of 
which no two are adjacent, is a well known NP-hard combinatorial optimization 
problem.
Since this problem has several applications in real life, there is need to find 
efficient algorithms to solve this problem.
Recently, Gaar and Rendl enhanced semidefinite programming approaches to tighten the upper 
bound given by the Lovász theta function. This is done by carefully selecting 
some 
so-called exact subgraph constraints (ESC)
and adding them to the semidefinite program of computing the Lovász theta 
function.

First, we provide two new relaxations that allow to compute the bounds faster 
without substantial 
loss of the quality of the bounds.
One of these two relaxations is based on including violated facets 
of the polytope representing the ESCs, 
the other one adds separating hyperplanes for that polytope.

Furthermore, we 
implement a branch and bound (B\&B) algorithm 
using these tightened relaxations in our bounding routine.
We compare the efficiency of our B\&B algorithm using the different 
upper bounds.
It turns out that already the bounds of Gaar and Rendl drastically reduce the 
number of nodes to be explored in the B\&B tree as compared to the Lovász theta 
bound.
However, this comes with a high computational cost.
Our new relaxations improve the run time of the overall B\&B algorithm, while 
keeping the number of nodes in the B\&B tree small.
\end{abstract}


\section{Introduction}

The stable set problem is a fundamental combinatorial optimization
problem. It is capable of modeling other combinatorial
problems as well as real-world applications and is therefore widely
applied in areas like operations research or computer science.
We refer to the survey~\cite{surveyWu} for more information and a review
of exact and heuristic algorithms. Most of the exact algorithms are
based on branch and bound (B\&B) and differ mainly by different upper and
lower bound computations.
\elli{ A recent paper using a MIP solver is
e.g.~\cite{Letchford2020}. The models used in that paper yield
computation times from less than a second up to half an hour on a
selection of DIMACS instances.

Outstanding results are obtained by an algorithm of Depolli
et. al~\cite{Depolli2013}. They introduced an algorithm using parallel
computing for finding maximum cliques in the context of protein
design. The algorithm consists of carefully implemented algorithmic
building blocks such as an approximate coloring algorithm, an initial
vertex ordering algorithm and the use of bit-strings for encoding the
adjacency matrix.}

In the 2015 survey~\cite{surveyWu}, no
exact algorithms using semidefinite programming (SDP) are mentioned.
One reason for the rare literature on SDP based B\&B
algorithms is the high computational cost for computing these bounds.
In this work we introduce an SDP based B\&B algorithm.
We formulate new SDP relaxations and develop solution algorithms to
compute these bounds with moderate computational expense, making them 
applicable within a B\&B scheme.

Before introducing the stable set problem, sometimes also referred to as 
vertex packing problem, we give the definition of a stable set.
Let $G = (V,E)$ be a simple undirected graph with $\abs{V} = n$ vertices 
    and $\abs{E} = m$ edges. A set $S \subseteq V$ is called stable if no 
    vertices in $S$ are adjacent.
    $S$ is called a maximal stable set if it is not possible to add a vertex to 
    $S$ without losing the stability property.
    The stability number $\alpha\left(G\right)$ denotes the maximum size of a 
    stable set in $G$, where size means the cardinality of the set.
    A stable set $S$ is called a maximum stable set if it has size $\alpha(G)$.

For convenience, from now on we always label the vertices of a graph with $n$ 
vertices
from~$1$ to~$n$. 
Computing $\alpha(G)$ can be done by solving the following optimization 
problem. 
\begin{alignat}{3}
\label{align:stableSetProblem}
\alpha(G) =
&\max & \sum_{i=1}^{n} x_i \\ \nonumber
&\st & x_i + x_j \leq 1  & ~\forall \set{i,j} \in E(G) \\ 
\nonumber
&& x \in \set{0,1}^n
\end{alignat}
    For a graph $G = (V,E)$,
the set of all 
stable set 
vectors  $\mathcal{S}(G)$ and the stable set polytope $\STAB(G)$ are 
defined as 
\begin{align*}
\mathcal{S}(G) &= \left\{ s \in \{0,1\}^{n} : s_{i}s_{j} = 0 \quad \forall 
\{i,j\} \in E\right\} \text{ and}\\ 
\STAB(G) &= \conv\left\{s : s \in \mathcal{S}(G)\right\}.
\end{align*}

Determining $\alpha(G)$ is NP-complete and the decision problem is among Karp's 
 21~NP-complete problems~\cite{karp}.
Furthermore, H{\aa}stad~\cite{hastad} proved
that $\alpha(G)$ is not approximable within $n^{1−\varepsilon}$ for any 
$\varepsilon > 0$ unless P=NP.
A well known upper bound on~$\alpha(G)$ is the Lov\'asz theta function 
$\vartheta(G)$.
Gr\"otschel, Lov\'asz and 
Schrijver~\cite{OurUsedFormOfLovasTheta} proved that
\begin{alignat}{3}
\label{lovaszTheta}
\vartheta(G) =
& \max & \allones{n}^{\transposedVec}x \phantom{iiii}  \\ \nonumber
& \st &  \diag(X) &= x\\ \nonumber
&& X_{i,j} &= 0 && \forall \{i,j\} \in E \\ \nonumber                  
&& \left(\begin{array}{cc}
1  & x^{\transposedVec} \\
x & X
\end{array}\right) &\succcurlyeq 0\\ \nonumber
&& X \in \Sym{n}&,~x \in \mathbb{R}^{n}
\end{alignat}
and hence provided a semidefinite program (SDP) to compute $\vartheta(G)$.
We define the feasible region of~\eqref{lovaszTheta} as
\begin{align*}
\ThetaBody^{2}(G) &= \left\{ (x,X) \in \R^{n}\times\Sym{n} \colon \diag(X) = x, 
\right.\\ 
& \left. \qquad \qquad X_{i,j}=0 \quad \forall \{i,j\} \in E, \quad
X-xx^{\transposedVec} 
\succcurlyeq 0 
\right\}.
\end{align*}
Clearly for each element $(x,X)$ of $\ThetaBody^{2}(G)$ the 
projection of $X$ onto its main diagonal is~$x$. The set of all projections 
\begin{equation*} 
\ThetaBody(G) = \left\{ x \in \R^{n}\colon \exists X \in \Sym{n} : (x,X) \in 
\ThetaBody^{2}(G) \right\}
\end{equation*}
is called theta body. More information on $\ThetaBody(G)$ can be found for 
example in Conforti, Cornuejols and 
Zambelli~\cite{ConfortiCornuejolsZambelliIntegerProgramming}.
It is easy to see that
$
\STAB(G) \subseteq \ThetaBody(G)
$
holds for every graph $G$, see~\cite{OurUsedFormOfLovasTheta}. 
Thus $\vartheta(G)$ is a relaxation of $\alpha(G)$.

This paper is structured as follows. In
Section~\ref{section:upperBoundESC} we introduce two new relaxations
using the concept of exact subgraph constraints. A branch and bound
algorithm that uses these relaxations is described in
Section~\ref{section:bab}, followed by the discussion of numerical
results in Section~\ref{sec:computationalExperiments}. 
Section~\ref{section:conclusions} concludes this paper.

\section{New Relaxations of the Exact Subgraph Constraints}
\label{section:upperBoundESC}

In this section we present two new approaches to find upper bounds 
on the stability number~$\alpha(G)$ of a graph $G$ starting from the Lovász 
theta 
function~$\vartheta(G)$ \elli{with so-called exact subgraph constraints, one 
based on 
violated facets and one based on 
separating hyperplanes. After introducing these approaches, we compare 
them both 
theoretically and practically. }

\subsection{Basic Setup for Exact Subgraph Constraints}
\label{sec:IEConvexCombination}
\elli{
Our approach is based on the idea of exact subgraph constraints that 
goes back to Adams, Anjos, Rendl and Wiegele~\cite{AARW} for combinatorial 
optimization problems that have an SDP relaxation and was recently 
computationally investigated by Gaar and 
Rendl~\cite{GaarRendlIPCO,GaarRendlFull} for the stable set, the Max-Cut and 
the coloring problem as a basis.
Starting from this, we present two relaxations of including  
exact subgraph constraints into the SDP for calculating~$\vartheta(G)$ that are 
computationally more 
efficient.
}

We first recapitulate the basic concepts of exact subgraph constraints with the 
notation 
from~\cite{GaarVersionsESH}. 
An upper bound on $\alpha(G)$ is given by the Lovász theta 
function~$\vartheta(G)$. Due to the SDP formulation~\eqref{lovaszTheta} 
it can be computed in polynomial time. 
Adams, Anjos, Rendl and Wiegele~\cite{AARW} proposed to 
improve~$\vartheta(G)$ as an upper bound 
by adding so-called exact subgraph constraints. These exact subgraph 
constraints can be used 
to strengthen SDP relaxations of combinatorial 
optimization problems with a certain property by including subgraph 
information. For the stable set problem we need the following definitions in 
order to introduce the exact subgraph constraints. 
For a graph $G$ the squared stable set polytope $\STAB^{2}(G)$ is defined as 
    \begin{equation*} 
\STAB^{2}(G) = \conv\left\{ss^{\transposedVec}: s \in 
\mathcal{S}(G)\right\}
\end{equation*}  
and matrices of the form $ss^{\transposedVec}$ for $s \in 
\mathcal{S}(G)$ are called stable set matrices.
Let $G_I$ denote the subgraph induced by the vertex set $I \subseteq V(G)$ with 
$|I| = \kI$.
With $X_I$ we denote the submatrix of~$X$ that results when we delete each row 
and column corresponding to a vertex that is not in $I$. In other words, $X_I$ 
is the submatrix of $X$ where we only choose the rows and columns corresponding 
to the vertices in $I$.
Then the constraint that asks the submatrix 
$X_I$ of~\eqref{lovaszTheta} for an induced subgraph $G_I$ to be in the squared 
stable set polytope 
$\text{STAB}^2(G_I)$ is called exact subgraph constraint (ESC).

The $k$-th level of the exact subgraph hierarchy introduced in~\cite{AARW} is 
the 
Lovász theta function~\eqref{lovaszTheta} with additional ESC for each 
subgraph of order $k$. In~\cite{GaarRendlIPCO,GaarRendlFull}  
this hierarchy is exploited computationally by including the ESC only for a set $J$  
of 
subgraphs and then considering
\begin{equation}
\label{SDPRelaxationWithESC}
\zESCJ = 
\max \left\{
\allones{n}^{\transposedVec}x:
(x,X) \in \ThetaBody^2(G),~
\XI \in \STAB^{2}(\GI) \quad \forall I \in J
\right\}. 
\end{equation}
Clearly, $\alpha(G) \leqslant \zESCJ$ holds for every set $J$ of subsets of 
$V(G)$, so $\zESCJ$ is an upper bound on $\alpha(G)$. 
One of the key remaining questions is how to solve~\eqref{SDPRelaxationWithESC}. 
We will compare different implementations and relaxations of this problem in 
the rest of the paper and start by considering existing methods.


The most straightforward way to solve~\eqref{SDPRelaxationWithESC} is 
to include the ESCs in a convex hull formulation as presented  
in~\cite{GaarRendlIPCO,GaarRendlFull,GaarVersionsESH}. 
We now recall the basic features and follow 
the 
presentation from~\cite{GaarVersionsESH}. 
As the ESC for a subgraph $\GI$ makes sure that  
$
\XI \in \STAB^{2}(\GI)
$
holds and the polytope $\STAB^{2}(\GI)$ is defined as the convex hull 
of the stable set matrices, the most 
intuitive way to formulate the ESC is as a convex combination. 
Towards that end, for a subgraph $\GI$ of $G$ induced by the 
subset $I \subseteq V$, let $|\mathcal{S}(\GI)| = \tI$ and let $\mathcal{S}(\GI) 
= 
\left\{\sIi{1}, \dots, \sIi{\tI}\right\}$.
Then the $i$-th stable set matrix $\SIi$ of the subgraph $\GI$ is 
defined 
as $\SIi = \sIi{i}(\sIi{i})^{\transposedVec}$.
As a result, the ESC $\XI \in \STAB^{2}(\GI)$ can be rewritten as
$$
\XI \in  \conv \left\{ \SIi: 
1 \leqslant i \leqslant \tI \right\}
$$
and it is natural to implement the ESC for the subgraph $\GI$ as
\begin{align*}
\XI = \sum_{i=1}^{\tI} 
\lambdaIi{i}\SIi, \quad \lambdaI \in \simplex_{\tI}, 
\end{align*}
where $\simplex_{\tI}$ is the $\tI$-dimensional simplex.

This means that when including the ESC for the subgraph $\GI$ 
into~\eqref{lovaszTheta} we have $\tI$ 
additional 
non-negative variables, one additional linear equality constraint for 
$\lambdaI$ and the matrix equality constraint which couples $\XI$ and 
$\lambdaI$.
We denote the number of equality constraints that are induced by the matrix 
equality constraint by $\bI$ and  note that $\bI \leqslant \binom{\kI+1}{2}$ 
holds.
With this formulation~\eqref{SDPRelaxationWithESC} can equivalently be written 
as
\begin{equation}
\label{SDPRelaxationWithESCConvHullOrig}
\ofvSDPescC = 
\max \left\{
\allones{n}^{\transposedVec}x:
(x,X) \in \ThetaBody^2(G),~
\XI = \sum_{i=1}^{\tI} \lambdaIi{i}\SIi, \quad \lambdaI \in \simplex_{\tI} 
\quad \forall I \in J
\right\},  
\end{equation}
so
$
\zESCJ = \ofvSDPescC
$
holds.
In practice, this SDP can be solved by interior point methods as long
as the number of ESC constraints is of moderate size.

Due to the fact that~\eqref{SDPRelaxationWithESCConvHullOrig} becomes a huge 
SDP as soon as the number of ESCs $|J|$ becomes 
large, Gaar and Rendl~\cite{GaarRendlIPCO,GaarRendlFull} proposed to use the 
bundle method to 
solve this SDP. The bundle method is an iterative procedure to find a global 
minimum of a non-smooth convex function and has been adapted for SDPs by 
Helmberg and Rendl~\cite{HelmbergRendlBundleSDP}.
As we use the bundle method only as a tool and do not enhance it any further, 
we refrain from presenting details here.

\subsection{Relaxation Based on Inequalities that Represent Violated Facets}
\label{sec:IEFacets}

We will see later on that the computational costs of a B\&B algorithm are 
enormous in the original version with 
the convex hull formulation~\eqref{SDPRelaxationWithESCConvHullOrig} and they 
are 
still substantial with the bundle approach  
from~\cite{GaarRendlIPCO,GaarRendlFull}. Therefore, we suggest two 
alternatives.

First, we present a relaxation of calculating the Lov\'{a}sz theta 
function with ESCs~\eqref{SDPRelaxationWithESC} that has already been mentioned 
in~\cite{GaarVersionsESH}, but has never 
been computationally exploited so far.
The key ingredient for this relaxation is the following observation. The 
polytope $\STAB^{2}(\GI)$ 
is given by its extreme points, which are the stable set matrices of $\GI$. Due 
to Weyl's theorem (see for 
example~\cite{NemhauserWolsey}) it can also be represented by its facets. This 
means that there are (finitely many) inequalities, such that the constraint 
$\XI \in \STAB^{2}(\GI)$ can be represented by these inequalities.

However, the facets and hence the inequalities depend on the stable set 
matrices and therefore on the subgraph $\GI$. Thus different subgraphs need 
different calculations that will 
lead to different inequalities. Gaar~\cite[\elli{Lemma 3}]{GaarVersionsESH} 
showed 
that 
adding
the ESC $ \XI \in \STAB^{2}(\GI)$
to the SDP calculating the 
Lov\'{a}sz theta 
function~\eqref{lovaszTheta} is equivalent
to adding the constraint  
$\XI \in \STAB^{2}(\GZeroi{\kI})$ where $\GZeroi{\kI} = (\VZeroi{\kI},\EZero)$ 
with $\VZeroi{\kI} = \{1, 
\dots, \kI\}$ and $\EZero = \emptyset$.

This implies that it is enough to calculate the facets of 
$\STAB^{2}(\GZeroi{\kI})$ and include these 
facets for each subgraph $\GI$ on $\kI$ vertices, instead of calculating the 
facets of $\STAB^{2}(\GI)$ for each subgraph $\GI$ separately. 
Let $\nrFacetsk{\elli{\kI}}$ be the number of facets of 
$\STAB^{2}(\GZeroi{\elli{\kI}})$ and let $\Fik{\elli{\kI}} \in \R^{\elli{\kI} 
	\times \elli{\kI}}$, $\fik{\elli{\kI}} \in \R$ 
for $1 \leqslant i \leqslant \nrFacetsk{\elli{\kI}}$ such that
\begin{align*}
\STAB^{2}(\GZeroi{\elli{\kI}}) = \left\{X \in \R^{\elli{\kI} \times 
	\elli{\kI}}: \left\langle \Fik{\elli{\kI}}, X 
\right\rangle \leqslant \fik{\elli{\kI}} \quad \forall 1 \leqslant i \leqslant 
\nrFacetsk{\elli{\kI}} \right\},
\end{align*}
so $(\Fik{\elli{\kI}},\fik{\elli{\kI}})$ is an inequality representing the 
$i$-th facet of 
$\STAB^{2}(\GZeroi{\elli{\kI}})$.
\elli{ 
We obtained $\nrFacetsk{\kI}$, $\Fik{\kI}$ and $\fik{\kI}$ for 
$\kI \leqslant 6$ in the way suggested in~\cite{GaarVersionsESH}.
For $\kI \geqslant 7$ this computation 
is beyond reach, as~$\nrFacetsk{7}$ is conjectured to be 
217093472~\cite{HPCutFacets}.
}

If we \elli{would} include all facets of $\STAB^{2}(\GZeroi{\kI})$ for each 
subgraph 
$\GI$ to replace the ESCs in~\eqref{SDPRelaxationWithESC}, then we would 
include a huge number of inequalities 
(\elli{$\nrFacetsk{5} = 368$ and 
$\nrFacetsk{6}= 116764$})
 and reach the limits of computing power 
rather soon. 
In order to reduce the number of inequalities, for each subgraph we 
 include only those inequalities that represent facets that are violated by 
 the current solution \elli{$\XOpt$}.
\elli{To be more precise, let~$\XOpt$ be the optimal solution 
of~\eqref{SDPRelaxationWithESC} for $J = \emptyset$,
i.e., the optimal solution of calculating the  
Lov\'{a}sz theta function.}
Then we 
define the indices of significantly violated 
facets of 
$\GI$\elli{, 
i.e., facets where the corresponding inequalities are violated at least by 
$\epsF$,}
 as
\begin{align*}
\indicesVFIPrime = \left\{1 \leqslant i \leqslant \nrFacetsk{\kI}:  
\left\langle \Fik{\kI}, \XIOpt \right\rangle > \fik{\kI} + \epsF\right\},
\end{align*}
where $\epsF$ is a small constant to take care of numerical inaccuracies of 
calculating $\XOpt$. 

Now we can further reduce the number of included inequalities in the following 
way. Although all $(\Fik{\kI},\fik{\kI})$ are different for different values of 
$i$, it could happen that for a subgraph $\GI$ there exist $1 \leqslant i \neq 
i' 
\leqslant \nrFacetsk{\kI}$ such that $(\Fik{\kI},\fik{\kI})$ and 
$(\FiPrimek{\kI},\fiPrimek{\kI})$ induce the same inequality. This is possible 
because they might differ only in positions $(j,j')$ with $j,j' \in I$ and 
$\{j,j'\} \in E$. Therefore, these different entries are multiplied with zero 
due to $[\XIOpt]_{j,j'} = 0$. Hence, let $\indicesVFI \subseteq 
\indicesVFIPrime$ be a set such that only one index among all indices in 
$\indicesVFIPrime$ which induce the same inequality is in~$\indicesVFI$.
Then we obtain the following relaxation of~\eqref{SDPRelaxationWithESC}, 
in which we 
include only inequalities that induce 
significantly violated facets of $\GI$
\begin{equation}
\label{SDPRelaxationWithESCFacetInequalities}
\ofvSDPescF = 
\max \left\{
\allones{n}^{\transposedVec}x:
(x,X) \in \ThetaBody^2(G),~
\left\langle \Fik{\kI}, \XI \right\rangle \leqslant \fik{\kI} \quad  
\forall i \in \indicesVFI \quad \forall I \in J
\right\}.
\end{equation}

Unfortunately for $\elli{\kI} \geqslant 7$ it is not possible to store and check the 
facets of $\STAB^{2}(\GZeroi{\elli{\kI}})$
for violation in reasonable memory and time due to  the huge number of facets.
Hence, we can perform this relaxation only for 
subgraphs $\GI$ of order $\kI \leqslant 6$.

%
%
%
%
%
%

\subsection{Relaxation Based on Separating Hyperplanes}
\label{sec:IESeparatingHyperplanes}
Next we consider  
another approach to implement a relaxation 
of~\eqref{SDPRelaxationWithESC} which can also be used for 
subgraphs $\GI$ of order $\kI \geqslant 7$ and which is based on including 
separating 
hyperplanes.

It uses 
the following fact.
Let \elli{$\Xtilde$ be any matrix in $\in \Sym{n}$ and 
let} $\projXI$ 
be the projection of \elli{$\XtildeI$} onto $\STAB^{2}(\GI)$. Then we can calculate the 
projection distance \elli{of $\Xtilde$ to $\STAB^{2}(\GI)$} as 
\begin{align}
\norm{\projXI - \elli{\XtildeI}}_{\elli{F}}^2\nonumber
&= \min_{\lambdaI \in \simplex_{\tI}} \norm{ \left(\sum_{i = 1}^{\tI} 
\lambdaIi{i}\SIi\right) - \elli{\XtildeI}}_{\elli{F}}^2\nonumber
= \min_{\lambdaI \in \simplex_{\tI}} \norm{ \sum_{i = 1}^{\tI} 
\lambdaIi{i}(\SIi - \elli{\XtildeI})}^{2}_{\elli{F}}\nonumber\\
&= \min_{\lambdaI \in \simplex_{\tI}} 
\sum_{j = 1}^{\kI} \sum_{j' = 1}^{\kI} 
\left(\sum_{i = 1}^{\tI} \lambdaIi{i} \left[  
\SIi - \elli{\XtildeI}
\right]_{j,j'}\right)^{2}\nonumber\\
&= \min_{\lambdaI \in \simplex_{\tI}} 
\sum_{j = 1}^{\kI} \sum_{j' = 1}^{\kI} 
\left(
\sum_{i = 1}^{\tI} \sum_{i'=1}^{\tI} \lambdaIi{i}\lambdaIi{i'} \left[  
\SIi - \elli{\XtildeI}
\right]_{j,j'}  \left[  
\SIiPrime - \elli{\XtildeI}
\right]_{j,j'}
\right)\nonumber\\
&= \min_{\lambdaI \in \simplex_{\tI}} 
\sum_{i = 1}^{\tI} \sum_{i'=1}^{\tI} \lambdaIi{i}\lambdaIi{i'}
\left(
\sum_{j = 1}^{\kI} \sum_{j' = 1}^{\kI} 
\left[  
\SIi - \elli{\XtildeI}
\right]_{j,j'}  \left[  
\SIiPrime - \elli{\XtildeI}
\right]_{j,j'}
\right)\nonumber\\    
&= \min_{\lambdaI \in \simplex_{\tI}} 
\lambdaI^{\transposedVec} \QI \lambdaI, \label{problemProjectionDistance}
\end{align}
where $\QI \in \R^{\tI \times \tI}$ and $[\QI]_{i,i'} = \left\langle \SIi - 
\elli{\XtildeI}, \SIiPrime - \elli{\XtildeI} \right\rangle$. $\QI$ is symmetric and positive 
semidefinite because it is a Gram matrix, so~\eqref{problemProjectionDistance} 
is a convex-quadratic program 
with $\tI$ variables, a convex-quadratic objective function and one linear 
equality 
constraint. With the 
optimal solution $\lambdaI$ of~\eqref{problemProjectionDistance} the projection 
of $\elli{\XtildeI}$ onto 
$\STAB^{2}(\GI)$ can be obtained by $\projXI = \sum_{i = 1}^{\tI} 
\lambdaIi{i}\SIi$. 
By defining
\begin{align*}
\HI = \frac{1}{\norm{\elli{\XtildeI}- \projXI}_{\elli{F}}} \left(\elli{\XtildeI}- \projXI \right) \quad 
\text{ and } \quad \hI = \frac{1}{\norm{\elli{\XtildeI}- \projXI}_{\elli{F}}} \left\langle 
\elli{\XtildeI} - \projXI, \projXI \right\rangle
\end{align*}
due to the separating hyperplane theorem (see for example Boyd and 
Vandenberghe~\cite{BoydVandenbergheConvexOptimization}) 
\begin{align}
\label{escSHForumlation}
\left\langle \HI, \XI \right\rangle \leqslant \hI
\end{align}
is a hyperplane that separates $\elli{\XtildeI}$ from $\STAB^{2}(\GI)$ such 
that $\XI = 
\projXI$ fulfills the inequality with equality.
Obviously~\eqref{escSHForumlation} is a relaxation of the 
ESC $\XI \in \STAB^{2}(\GI)$, so
\begin{equation}
\label{SDPRelaxationWithESCSeparatingHyperplane}
\ofvSDPescH = 
\max \left\{
\allones{n}^{\transposedVec}x:
(x,X) \in \ThetaBody^2(G),~
\left\langle \HI, \XI \right\rangle \leqslant \hI \quad \forall I \in J
\right\}
\end{equation}
is another relaxation of~\eqref{SDPRelaxationWithESC}
\elli{that depends on the chosen $\Xtilde$}.

\subsection{Theoretical Comparison  of the Relaxations}
\label{sec:IEComparison}
We briefly comment on some theoretical properties of  
$\ofvSDPescC$, $\ofvSDPescF$ and $\ofvSDPescH$.
We start by analyzing the upper bounds we obtain. Due to the fact that 
$\ofvSDPescF$ and $\ofvSDPescH$ are relaxations of $\ofvSDPescC$, we know that
\begin{align*}
\alpha(G) \leqslant \ofvSDPescC \leqslant \ofvSDPescF, \ofvSDPescH \leqslant 
\vartheta(G)
\end{align*}
holds for every graph $G$ and every set $J$. 

Another important observation is the following. Whenever we include the 
ESC of the subgraph~$\GI$ into the SDP \elli{ computing $\ofvSDPescC$,
the stable set problem is solved exactly on this subgraph~$\GI$.  
However, when computing
$\ofvSDPescF$ and $\ofvSDPescH$ we do not include the ESC but only a
relaxed version of it.
Hence, in the optimal solutions of these two
relaxations, 
it could still be the case that the ESC is not fulfilled, i.e.,  
for the subgraph $\GI$ we do not have an exact solution.
Hence, it is possible that we still find violated inequalities (representing facets or hyperplanes) in these cases. As a consequence, for $\ofvSDPescC$ it does not 
make sense to include the ESC for the same 
subgraph twice, but for $\ofvSDPescF$ and $\ofvSDPescH$ it is 
possible 
that we want to include a relaxation of the very same ESC twice} with different facets 
or a different separating hyperplane.

Finally let us consider the sizes of the SDPs to solve. In all three versions 
$\ofvSDPescC$, $\ofvSDPescF$ and $\ofvSDPescH$ we solve the SDP 
of the Lov\'{a}sz theta function~\eqref{lovaszTheta} with additional 
constraints, so in all three SDPs we have a matrix variable of dimension $n+1$ 
which has to be positive semidefinite \melanie{(psd)} and $n+m+1$ linear 
equality constraints. 
Additionally to that we have  $\sum_{I \in J} \tI$ non-negative variables
and $\nrESC+\sum_{I \in J}\bI$ equality constraints for~$\ofvSDPescC$, 
$\sum_{I \in J} |\indicesVFI|$ inequalities for $\ofvSDPescF$, 
and $\nrESC$ inequalities for $\ofvSDPescH$.
\melanie{Table~\ref{tab:sizeSDPs} gives an overview of the different sizes of the SDPs.}

\begin{table}[htbp]
	\caption{\melanie{Sizes of the SDPs to compute $\ofvSDPescC$, 
		$\ofvSDPescF$ and $\ofvSDPescH$}}
	\label{tab:sizeSDPs}
	\setlength{\tabcolsep}{5pt}    
	\begin{center}
		\begin{tabular}{l|ccc}
			& \multicolumn{1}{c}{$\ofvSDPescC$} 
			& \multicolumn{1}{c}{$\ofvSDPescF$}
			& \multicolumn{1}{c}{$\ofvSDPescH$} 
			\\  \\[-0.75em] \hline  \\[-0.75em]
			\melanie{dimension psd matrix variable}
			&$n+1$
			&$n+1$
			&$n+1$
			\\  \\[-0.75em]
			\melanie{\# non-negative variables}
			& $\sum_{I \in J} \tI$
			& 0
			& 0
			\\ 
			\\[-0.75em]
			\melanie{\# linear equality constraints}
			&$n+m+1 + \nrESC+\sum_{I \in J}\bI$ 
			&$n+m+1$ 
			&$n+m+1$ 
			\\ 
			\\[-0.75em]
			\melanie{\# linear inequality constraints}
			&0
			& $\sum_{I \in J} |\indicesVFI|$
			& $\nrESC$
		\end{tabular}
	\end{center}
\end{table}

\subsection{Computational Comparison of the Relaxations}
\label{sec:SmallComputationalComparisonVFSH}
Before we perform a large scale comparison 
of $\ofvSDPescC$, $\ofvSDPescF$ and $\ofvSDPescH$
within a B\&B algorithm, 
we investigate a small example.
\begin{example}
    \label{Ex:OneGraph}
    We consider a random graph $G=G_{100,15}$ from the Erd\H{o}s-Rényi model 
    $G(n,p)$ with $n=100$ and $p=0.15$. A random graph from this model has $n$ 
    vertices and every edge is present with probability $p$ 
    independently from all other edges.
    For the chosen graph, $\vartheta(G_{100,15}) = 27.2003$ and 
    $\alpha(G_{100,15}) = 24$ holds, so 
    $$
    24 \leqslant \ofvSDPescC \leqslant 
    \ofvSDPescF, \ofvSDPescH \leqslant  27.2003
    $$ 
    holds for every set $J$.

    All the computations were performed on an Intel(R) Core(TM) 
    i7-7700 CPU @ 3.60GHz with 32 GB RAM with the MATLAB version R2016b and 
    with 
    MOSEK version 8.
    \elli{
    In the computations, we use $\epsF = 0.00005$ 
	and we include a separating hyperplane for a subgraph whenever 
	the projection distance is greater or equal to $0.00005$.    
	}
    
    We compute $\ofvSDPescC$, $\ofvSDPescF$ and $\ofvSDPescH$ for 
    different sets $J$, which all consist of subsets of vertices 
    of 
    size~$\kI = 5$ and only differ in the number \elli{$\nrESCNew=\nrESC$} of 
    included 
    ESCs.
    \elli{To be more precise, we consider five different sets 
    $J = \elli{\JbNew}$ 
    with $\nrESCNew = |\JbNew| \in \{221,443,664,886,1107\}$.
	These values of $\nrESCNew$ are chosen in such a way that
	the number of linear 
	equality constraints which are induced by the 
	matrix equalities from the ESCs in the convex hull formulation, 
	i.e., $\sum_{I \in \JbNew} \bI$, is in $\{3000,6000,9000,12000,15000\}$.
	To choose the subsets in $\JbNew$, 
	we first determine~$\XOpt$ as the optimal solution 
	of~\eqref{SDPRelaxationWithESC} for $J = \emptyset$,
	i.e., the optimal solution of calculating the  
	Lov\'{a}sz theta function~\eqref{lovaszTheta}.
	Then we generated~$3\nrESCNew$ subgraphs $\GI$ of order $\kI$ randomly and 
	included those $\nrESCNew$ subsets $I$ into $\JbNew$, where the corresponding $\XOptI$ have 
	the largest projection distances to $\STAB^{2}(\GI)$.
	}
	\elli{For computing $\ofvSDPescH$ we choose $\Xtilde = \XOpt$.}
    

\begin{table}[ht]
    \caption{The values of $\ofvSDPescC$, $\ofvSDPescF$ and 
    $\ofvSDPescH$ for   
    $G=G_{100,15}$ for different sets $\elli{\JbNew}$
}
\label{tab:tightenedThetaFunctionExample}    
    \begin{center}
        \begin{tabular}{l|rrrrr}
            & \multicolumn{1}{c}{$\elli{J_{221}}$} 
            & \multicolumn{1}{c}{$\elli{J_{443}}$} 
            & \multicolumn{1}{c}{$\elli{J_{664}}$} 
            & \multicolumn{1}{c}{$\elli{J_{886}}$} 
            & \multicolumn{1}{c}{$\elli{J_{1107}}$}
            \\ \\[-0.75em] \hline \\[-0.75em]
            $\ofvSDPescC$
            &26.9905
            &26.9299
            &26.8684 
            &26.8496
            &26.8278
            \\ 
            \\[-0.75em]
            $\ofvSDPescF$
            &26.9975
            &26.9393
            &26.8807
            &26.8602
            &26.8397
            \\ 
             \\[-0.75em]
            $\ofvSDPescH$
            &27.0104
            &26.9741
            &26.9215
            &26.8992
            &26.8898
            \\ 
        \end{tabular}
    \end{center}
\end{table}

If we consider Table~\ref{tab:tightenedThetaFunctionExample} with the improved 
upper bounds then we see that if 
$\elli{\nrESCNew}$ increases, all upper bounds  
$\ofvSDPescC$,  $\ofvSDPescF$ and $\ofvSDPescH$ improve.
Furthermore, one can observe that for a fixed set $\elli{\JbNew}$ the obtained 
bounds
of $\ofvSDPescC$ are best, those of $\ofvSDPescF$ are a little bit worse and 
those of $\ofvSDPescH$ are even a little bit more worse, 
\elli{i.e., empirically the bounds obtained by using $\ofvSDPescF$ are better 
than those coming from $\ofvSDPescH$ in our example.}

\begin{table}[htbp]
    \caption{The running times in seconds for computing $\ofvSDPescC$, 
    $\ofvSDPescF$ and $\ofvSDPescH$ of 
    Table~\ref{tab:tightenedThetaFunctionExample}
}
\label{tab:timeExample}    
    \begin{center}
        \begin{tabular}{l|rrrrr}
            & \multicolumn{1}{c}{$\elli{J_{221}}$} 
            & \multicolumn{1}{c}{$\elli{J_{443}}$} 
            & \multicolumn{1}{c}{$\elli{J_{664}}$} 
            & \multicolumn{1}{c}{$\elli{J_{886}}$} 
            & \multicolumn{1}{c}{$\elli{J_{1107}}$}
            \\  \\[-0.75em] \hline  \\[-0.75em]
            $\ofvSDPescC$
            &8.14
            &29.30
            &66.95
            &139.64
            &279.11
            \\ 
            \\[-0.75em]
            $\ofvSDPescF$
            &0.61 
            &1.15 
            &1.94 
            &2.72 
            &3.81
            \\ 
             \\[-0.75em]
            $\ofvSDPescH$
            &0.75
            &1.25
            &1.93 
            &2.57
            &3.34
            \\ 
        \end{tabular}
    \end{center}
\end{table}

    Next we consider the running times for computing $\ofvSDPescC$, 
        $\ofvSDPescF$ and $\ofvSDPescH$ in Table~\ref{tab:timeExample}. 
Here we see that the time it takes so solve 
$\ofvSDPescC$ is extremely high and increases drastically if the number of 
included ESCs gets larger. Both our relaxations $\ofvSDPescF$ and $\ofvSDPescH$ 
reduce the running times 
enormously. The running times for $\ofvSDPescF$ and $\ofvSDPescH$ are 
comparable, but computing $\ofvSDPescH$ is slightly faster for including a 
large number 
of ESCs as we only 
include one additional inequality in $\ofvSDPescH$ whereas we 
may include several inequalities that represent facets in $\ofvSDPescF$.    
    
\begin{table}[ht]
    \caption{The average projection distances of $\XI$ to 
    $\STAB^{2}(\GI)$ over all 
    $I \in \elli{\JbNew}$  before (i.e., \elli{$X = \XOpt$} is the optimal 
    solution of~\eqref{lovaszTheta}), and after 
    (i.e., \elli{$X \in \{\XOptC, \XOptF, \XOptH\}$} is the optimal solution of 
    $\ofvSDPescC$, $\ofvSDPescF$ and  
    $\ofvSDPescH$)
    including the ESCs
}
\label{tab:projectionDistanceExample}        
    \begin{center}
        \begin{tabular}{l|rrrrr}
            & \multicolumn{1}{c}{$\elli{J_{221}}$} 
            & \multicolumn{1}{c}{$\elli{J_{443}}$} 
            & \multicolumn{1}{c}{$\elli{J_{664}}$} 
            & \multicolumn{1}{c}{$\elli{J_{886}}$} 
            & \multicolumn{1}{c}{$\elli{J_{1107}}$}
            \\  \\[-0.75em] \hline  \\[-0.75em]
            before including ESCs 
            &0.03095
            &0.03014
            &0.03057
            &0.02982 
            &0.03032
            \\  \\[-1em] \hline  \\[-0.75em] 
            after computing $\ofvSDPescC$
            &0.00005
            &0.00004
            &0.00004
            &0.00004
            &0.00004
            \\ 
            \\[-0.75em]
            after computing $\ofvSDPescF$
            &0.00151 
            &0.00087 
            &0.00115
            &0.00080
            &0.00051
            \\ 
            \\[-0.75em]
            after computing $\ofvSDPescH$
             &0.00290
             &0.00256
             &0.00252
             &0.00196
             &0.00183
             \\ 
        \end{tabular}
    \end{center}
\end{table}

As a next step we investigate the projection 
distances.    
\elli{Recall that $\XOpt$ is the optimal solution of calculating the 
Lov\'{a}sz theta function~\eqref{lovaszTheta}}.
Let $\XOptC$, $\XOptF$ and $\XOptH$ be the optimal solution of 
calculating 
$\ofvSDPescC$, $\ofvSDPescF$ and  $\ofvSDPescH$,  
respectively.
In Table~\ref{tab:projectionDistanceExample} we see that the average projection 
distance of $\XOpt$ is significantly larger than $0$ before 
including the ESCs, so there are several violated ESCs. 
As soon as the ESCs are included the average projection distance for  
$\XOptC$ is almost zero, so the ESCs are almost satisfied. In theory 
they should all be zero, but as MOSEK is not an exact solver, the optimal 
solution is subject to numerical inaccuracies.
If we turn to $\ofvSDPescF$, then the projection distances of $\XOptFI$ are not 
as close to zero as those for $\XOptC$, because $\ofvSDPescF$ is only a 
relaxation of $\ofvSDPescC$.
Also the average projection distance of 
$\XOptHI$ after 
solving~$\ofvSDPescH$ is greater than 
the one obtained with $\ofvSDPescF$. This is in tune with the fact that 
the upper bounds obtained in the latter case are better \elli{for this 
instance}.
Furthermore, note that the 
average projection distances for $\XOptFI$ and $\XOptHI$
decrease as $\elli{\nrESCNew}$ increases. This is not 
surprising, as more ESCs mean that a bigger portion of the graph is 
forced into the stable set structure.

\begin{table}[ht]
    \caption{The average number of violated facets  $|\indicesVFI|$ over all 
    subgraphs $\GI$ 
    with $I \in \elli{\JbNew}$ before and after including the ESCs  
    and computing $\ofvSDPescF$}
\label{tab:averageNrOfViolatedFacetsExample}    
    \begin{center}
        \begin{tabular}{l|rrrrr}
            & \multicolumn{1}{c}{$\elli{J_{221}}$} 
            & \multicolumn{1}{c}{$\elli{J_{443}}$} 
            & \multicolumn{1}{c}{$\elli{J_{664}}$} 
            & \multicolumn{1}{c}{$\elli{J_{886}}$} 
            & \multicolumn{1}{c}{$\elli{J_{1107}}$}
            \\  \\[-0.75em] \hline  \\[-0.75em]
            before including ESCs 
            &1.53 
            &1.56
            &1.57
            &1.56
            &1.56
            \\  \\[-1em] \hline  \\[-0.75em]
            after computing $\ofvSDPescF$
            & 0.23
            & 0.14
            & 0.16
            & 0.15
            & 0.09
            \\  
        \end{tabular}
    \end{center}
\end{table}

Finally we present in 
Table~\ref{tab:averageNrOfViolatedFacetsExample} the average number of 
violated facets.
As one can see the 
average number of violated facets before including the ESCs is already very 
low. 
This means that we do not include too many inequalities that 
represent facets in the computation of~$\ofvSDPescF$.
Furthermore, 
the average number of facets that are violated by $\XOptF$ decreases 
significantly compared to the average number of violated facets before 
including the relaxations of the ESCs.
This is very encouraging because one could imagine a scenario where we 
iteratively add violated facets of one subgraph and then the optimal 
solution violates different facets. However, the computations suggest that this 
does not happen too often.
Like before in Table~\ref{tab:projectionDistanceExample} we see that the more 
ESCs are included, the more stable set structure is captured and therefore 
the less facets are violated after including the relaxations of the ESCs.
\hfill $\bigcirc$
\end{example}

Let us briefly summarize the key points of Example~\ref{Ex:OneGraph}.
Usually the upper bounds obtained by $\ofvSDPescF$ are only slightly worse than 
those of $\ofvSDPescC$, but the running times are only a fraction. 
Unfortunately, this approach works only for subgraphs of order at most~$6$. 
Also $\ofvSDPescH$ yields good upper bounds in slightly better running time 
than $\ofvSDPescF$, but the bounds are a little bit worse. A major benefit of 
this approach is that it can be used for subgraphs of any order.

In a nutshell, both $\ofvSDPescF$ and $\ofvSDPescH$ are relaxations of 
$\ofvSDPescC$ that reduce the running times drastically by worsening the bounds 
only a little bit. As a result these bounds are very promising for including 
them into a B\&B algorithm for stable set.

\section{Branch and Bound for the Stable Set Problem}
\label{section:bab}

The aim of this section is to present our implementation of an exact branch and 
bound (B\&B) algorithm for the stable set 
problem~\eqref{align:stableSetProblem}. 

\subsection{Our Branch and Bound Algorithm}
\label{sec:babalgo}
We start by detailing the general setup of our B\&B algorithm. 
Towards this end, keep in mind that in a solution 
of~\eqref{align:stableSetProblem} the binary 
variable $x_i$ is equal to $1$ if  vertex  $i$ is in the stable set, and 
$x_i = 0$ otherwise. 



For our B\&B algorithm for the stable set problem we choose a vertex $i \in 
V(G)$ and divide the  
optimization problem in 
a node  of the B\&B tree
into the subproblem where vertex $i$ is in the stable set (i.e., set the branching variable $x_i = 1$) and 
a 
second 
subproblem where $i$ is not in the stable set (i.e., set the branching variable 
$x_i = 0$).

It turns out that in each node of the B\&B tree the 
optimization problem is of the form
\begin{align}
\label{align:problem-in-node}
P(G,c) = 
c + \max \quad &  \sum\limits_{i \in V(G)} x_i   \\ \nonumber
\st \quad & x_i + x_j \leq 1 \quad  \forall \set{i,j} \in E(G)\\ 
\nonumber
& x_i \in \set{0,1} \quad \ \,\forall i \in V(G)
\end{align}
for some graph $G$, so in each node we have to solve a stable set problem and add 
a constant term $c$ to the objective function value. 
Indeed, \elli{by fixing a branching variable $x_i$ to $0$ or $1$,
we shrink the graph and create subproblems that are again
stable set problems of the  
form~\eqref{align:problem-in-node} but with a smaller graph and some offset $c$.}
To be more precise, for the 
subproblem with $x_i = 1$ the objective function value 
of~\eqref{align:stableSetProblem} increases by~$1$ because there is one more 
vertex in the stable set. 
Furthermore, 
all neighbors of $i$ can not be in the maximum stable set because $i$ is 
already in the stable set.
So we can set 
$x_j = 0$ for all $j \in N_G(i)$ if 
 $N_G(i) = \set{j 
    \in V(G) \mid \set{i,j} \in E(G) }$ 
denotes the set of neighbors of the vertex $i$ in~$G$.
Furthermore, we can delete $i$ and all neighbors of 
$i$ in the current graph $G$ and search for a maximum stable set in the new 
graph $G'$
of 
smaller order, where $G'= G[U']$ is the subgraph of~$G$ induced 
by $U' = V(G)\setminus \bigl(\set{i} \cup N_G(i) \bigr)$. 
Hence, the 
subproblem to solve in the new node has the form 
$P(G',c+1)$.

In the subproblem for $x_i = 0$ the vertex $i$ is not in 
the 
stable 
set. We can 
remove the vertex $i$ from the graph and search for a maximum stable set in the 
induced subgraph $G'' = G[U'']$ with $U'' = V(G)\setminus \set{i}$. This 
boils 
down 
to solving $P(G'',c)$ in the new node 
of the B\&B tree.

Note that whenever we delete a vertex $i$ from the graph in the branching 
process, we set the according variable $x_i$ to a fixed value. 
Consequently, in every node, all vertices of the original graph $G$ are either 
still present, or the value of the variable corresponding to them is implicitly 
fixed.
Furthermore, we exclude all non-feasible solutions by deleting all neighbors in 
case of setting the branching variable to~$1$.
Hence, every time we set a variable of a vertex to $1$ the set of all vertices 
of which the variable is set to $1$ remains stable and we only 
obtain feasible solutions of~\eqref{align:stableSetProblem}.
Therefore, from a feasible solution 
of~\eqref{align:problem-in-node} in any node we can determine a feasible 
solution of~\eqref{align:stableSetProblem} 
with the same objective 
function value.

The order of the graph to consider in a node shrinks whenever we 
branch.
As a consequence the 
B\&B tree is of finite size.
Whenever we reach a node with a suproblem on a graph with 
less or equal to 23 vertices, 
we 
solve the problem by \elli{
a fast enumeration procedure that can be employed whenever
the subproblems become sufficiently small.
To do so, we iterate  
easily -- and especially fast -- over} all subsets of $V$ in 
descending order 
with an implementation of Hinnant~\cite{combAndPerm} in \CPP.

\paragraph{Bounds}
We do not want to solve the 
subproblems~\eqref{align:problem-in-node} in each node to optimality, but only 
compute an   
upper and a lower bound on the optimal objective function value. This boils 
down to obtaining bounds on the stability number of the graph considered 
in~\eqref{align:problem-in-node}. We present details on lower bounds obtained 
by heuristics in Section~\ref{sec:lower-bound}.

As upper bounds we use the relaxations based on ESCs 
in four different versions, namely the convex hull formulation or the 
bundle method as detailed in Section~\ref{sec:IEConvexCombination}, 
the 
violated facets version as described in Section~\ref{sec:IEFacets}  
or 
the separating hyperplanes version as presented in 
Section~\ref{sec:IESeparatingHyperplanes}.  
For choosing the subset $J$ of ESCs, 
we follow the approach of~\cite{GaarRendlFull} in our computations and perform 
several cycles\elli{, i.e., iterations of the repeat until loop,} of solving 
\elli{(a relaxation) of}~\eqref{SDPRelaxationWithESC} and then adjusting 
the set $J$\elli{, as illustrated in Algorithm~\ref{alg:bounds}}.

\elli{In particular,} we start with $J = \emptyset$\elli{, as preliminary computations have shown that carrying over ESC to subproblems does not pay of,} and in each \elli{cycle we update $J$ depending on the current optimal solution~$\XHyp$ of the SDP solved.}
\elli{We remove all previously added ESCs where the associated dual variables of the optimal solution  have absolute value less than 0.01.}
For finding violated subgraphs (i.e., subgraphs for 
which the ESC does not hold \elli{in $\XHyp$}) we use the methods presented 
in~\cite{GaarRendlFull}, so we use a local search heuristic to find submatrices 
of $\XHyp$ that minimize the inner product with some matrices. We let the local 
search heuristic run for $9n$ times \melanie{and add random subgraphs to obtain $9n$ subgraphs without duplicates. From these subgraphs we add the $3n$ most violated ones (subsets $I$ with largest projection distance of $\XHyp_I$ to $\STAB^{2}(\GI)$)  to $J$.
}

To reduce computational effort, we stop cycling as soon as we do not expect 
to be able to prune within the next 5 cycles \elli{assuming that the decrease
of the upper bound $z^*$ in each future cycle is 0.75 of the average decrease we 
had in the previous cycles.}

\begin{algorithm}
\DontPrintSemicolon
\SetAlgoLined
\KwIn{Graph $G$ at current node, method {\em convex-hull}, {\em bundle-method}, {\em violated-factes}, or {\em separating-hyperplanes}}
\KwOut{Upper bound $z^*$}
$J=\emptyset$\;
\Repeat{optimistic forecast does not suggest pruning}{
  \lIf{convex-hull}{let $\XHyp$ be an optimal solution of \eqref{SDPRelaxationWithESCConvHullOrig} with objective function value $z^*$}
  \lIf{bundle-method}{let $\XHyp$ be the solution of approximately solving \eqref{SDPRelaxationWithESCConvHullOrig} using the bundle method with objective function value $z^*$}
  \lIf{violated-facets}{let $\XHyp$ be an optimal solution of \eqref{SDPRelaxationWithESCFacetInequalities}  with objective function value $z^*$}
  \lIf{separating-hyperplanes}{let $\XHyp$ be an optimal solution of 
  \eqref{SDPRelaxationWithESCSeparatingHyperplane} with objective function 
  value $z^*$, where we choose $\Xtilde$ as $\XHyp$ of the last cycle}
  remove ESCs from $J$ if associated dual variables have small absolute value\;
  search for ESCs violated by $\XHyp$ and add them to $J$\;}
\Return{$z^*$}
\caption{Upper bound computation at a B\&B node}
\label{alg:bounds}
\end{algorithm}

For our computations we use the implementation of the bundle 
method as it is 
detailed in~\cite{GaarRendlFull}, in particular with all 
specifications given 
in Section~6.3 therein and we let the bundle method run for at most 15 
iterations in each cycle.

\paragraph{Branching Rule}
An important question in the B\&B algorithm is how to choose the branching 
variable.
In our implementation we follow the approach to first deal with vertices, for 
which we know least whether they will be in a maximum stable set or not 
(\emph{``difficult first''}) in 
order to find an optimal solution soon. 

All our upper bounds are based on the Lov\'asz theta 
function~\eqref{lovaszTheta}, so we can use the intuition
that the closer an entry $x_i$ 
is to $1$ in an optimal solution of~\eqref{lovaszTheta}, the more likely it is 
that this vertex $i$ is in a maximum stable set.
Hence, we choose the variable $x_i$ with value
closest to $0.5$ as
branching variable. 

More on branching rules for the stable set problem can for 
example be found in~\cite{branching1,branching2}.

\paragraph{Diving Strategy}
We implemented a best first search 
strategy, 
where we always consider the open subproblem with the largest upper bound next.
We expect that we find a large stable set in this branch of 
the B\&B tree because the difference between the global lower bound  
and the upper bound for this branch is the highest of all.

\subsection{Heuristics to Find Large Stable Sets}
\label{sec:lower-bound}
It is  crucial to find a good lower bound on $\alpha(G)$ early in the B\&B 
algorithm to prevent the growth of the  B\&B  tree and 
therefore solve the stable set problem more efficiently.
In~\cite{krishna} and~\cite{surveyWu}, for example, one can find references to 
some heuristics to find a 
large stable set.
In our implementation we use several different heuristics. 

The first heuristic makes use of the vector $x$ from the SDP 
formulation~\eqref{lovaszTheta} of $\vartheta(G)$, which is available from 
the upper bound computation.
This vector consists of $n$ elements between $0$ and $1$. The value $x_i$ gives 
us some intuition about the $i$-th vertex of the graph, namely the closer it is 
to $1$, the more likely it is that the vertex is in a maximum stable set. 
Hence, we sort the vertices in descending order according to their value in $x$ 
and then 
add the vertices in this order to a set $S$, such that the vertices of $S$ 
always remain a stable set.  
In the following we refer to this heuristic as \mytag{HT}.

Furthermore, we use a heuristic introduced by 
Kahn I., Ahmad and Kahn M. in~\cite{hvs1} based on vertex covers and  
vertex supports.
A subset $C$ of the vertices of a graph is called vertex cover if for each edge 
at least one of the two incident vertices is in $C$. 
The vertex support of a vertex is defined as the sum of the degrees 
of all vertices in the neighborhood of this vertex. 
If $C$ is a vertex cover, then clearly $V(G)\setminus C$ is a stable set, so 
instead of searching for a maximum stable set we can search for a vertex cover 
of minimum cardinality. 
In a nutshell, the heuristic of~\cite{hvs1} searches for a vertex with maximum 
vertex support 
in the 
neighborhood of the vertices with minimum vertex support.
If there is more than one vertex with maximum support, one with maximum degree 
is chosen.
This vertex is added to the vertex cover $C$ and all incident edges are removed.
The above steps are repeated until there are no edges left in the graph.
In the end we obtain a hopefully large stable set with $V(G)\setminus C$.
We denote this heuristic by \mytag{HVC}.

Finally we use a heuristic proposed by Burer, Monteiro and Zhang 
in~\cite{burer02}.
Their heuristic is based on the SDP formulation of the Lovàsz theta function 
with additional restriction to the matrix variable to be of low rank.
With rank one, a local maximizer of the problem yields a maximal stable set, 
whereas with rank two the stable set corresponding to the local maximizer does 
not necessarily have to be maximal, but one can escape to a higher local 
maximizer which corresponds to a maximal stable set.
The C~source code of this heuristic is online available at~\cite{hvs1github}.
We use this code with the parameters rank set to $2$ and the number of 
so-called escapes set to $1$ in a first version and set to $5$ in a second 
version. 
Both parameter settings are among the choices that were tested 
in~\cite{burer02}.
We will refer to this versions with \mytag{H21} and \mytag{H25}.

In the B\&B algorithm we perform the heuristic \myref{H25} in the root node 
with a 
time limit of $20$ seconds. Then we only perform the heuristics in every third  
node of the B\&B tree. 
\myref{HT} and \myref{HVC} are very fast, so we let them run in each node we 
run heuristics.
Furthermore, in the first $10$ nodes  of the B\&B tree  we perform \myref{H25} 
with the 
hope 
to find a stable set of cardinality $\alpha(G)$ as soon as possible, but we do 
not allow the heuristic to iterate longer than $7$ seconds.
For graphs with less than $200$ vertices we additionally perform \myref{H21} 
with the running time limited to $1$ second.
On graphs with more vertices we perform with probability $0.05$ \myref{H25} 
and 
a time limit of $7$ seconds and in the other cases \myref{H21} with a time 
limit of $3$ seconds.
In a 
computational comparison in the master's thesis of 
Siebenhofer~\cite{melanie-diplomarbeit}, this turned out to be the best combination of heuristics.

\section{Computational Experiments}
\label{sec:computationalExperiments}

\begin{figure}[hbt]
	\centering
	\begin{subfigure}[b]{0.5\textwidth}
		\centering
			\begin{tikzpicture}[trim axis left] \footnotesize
			\begin{axis}[legend pos = north west, 
			xtick={0,2500,5000,7500,10000,12500,15000},
			ytick={0,0.2,0.4,0.6,0.8,1},
			enlarge x limits=false,
			enlarge y limits={upper},
			legend cell align={left},
			xlabel= time (sec),
			ylabel= \% of instances solved,
			ymin=0,
			ymax=1]
			\addplot+[red, mark=10-pointed star] table [x=time, y=CH]{performance-time.dat};
			\addlegendentry{\footnotesize CH}
			\addplot+[green, mark options={fill=green}] table [x=time, y=BD]{performance-time.dat};
			\addlegendentry{\footnotesize BD}
			\addplot+[blue, mark options={fill=blue}] table [x=time, y=VF]{performance-time.dat};
			\addlegendentry{\footnotesize VF}
			\addplot+[orange, mark options={fill=orange}] table [x=time, y=SH]{performance-time.dat};
			\addlegendentry{\footnotesize SH}
			\end{axis}
			\end{tikzpicture}
		\caption{Time to solve the instances}
	\end{subfigure}%
	\hfill
	\begin{subfigure}[b]{0.5\textwidth}
			\begin{tikzpicture} \footnotesize
			\begin{axis}[
			xtick={0,0.1,0.2,0.3,0.4,0.5,0.6},
			ytick={0,0.2,0.4,0.6,0.8,1},
			yticklabels=\empty,
			enlarge x limits=false,
			enlarge y limits={upper},
			xlabel= relative gap (ub-lb)/lb,
			ymin=0,
			ymax=1]
			\addplot+[red, mark=10-pointed star] table [x=gap, y=CH]{performance-gap-relative.dat};
			\addplot+[green, mark options={fill=green}] table [x=gap, y=BD]{performance-gap-relative.dat};
			\addplot+[blue, mark options={fill=blue}] table [x=gap, y=VF]{performance-gap-relative.dat};
			\addplot+[orange, mark options={fill=orange}] table [x=gap, y=SH]{performance-gap-relative.dat};
			\legend{};
			\end{axis}
			\end{tikzpicture}
		\caption{Gap after the time limit}
	\end{subfigure}
	\caption{The performance comparison of the bounds~\protect\myref{CH}, 
	\protect\myref{BD}, \protect\myref{VF} and~\protect\myref{SH}}
	\label{fig:performance}
\end{figure}
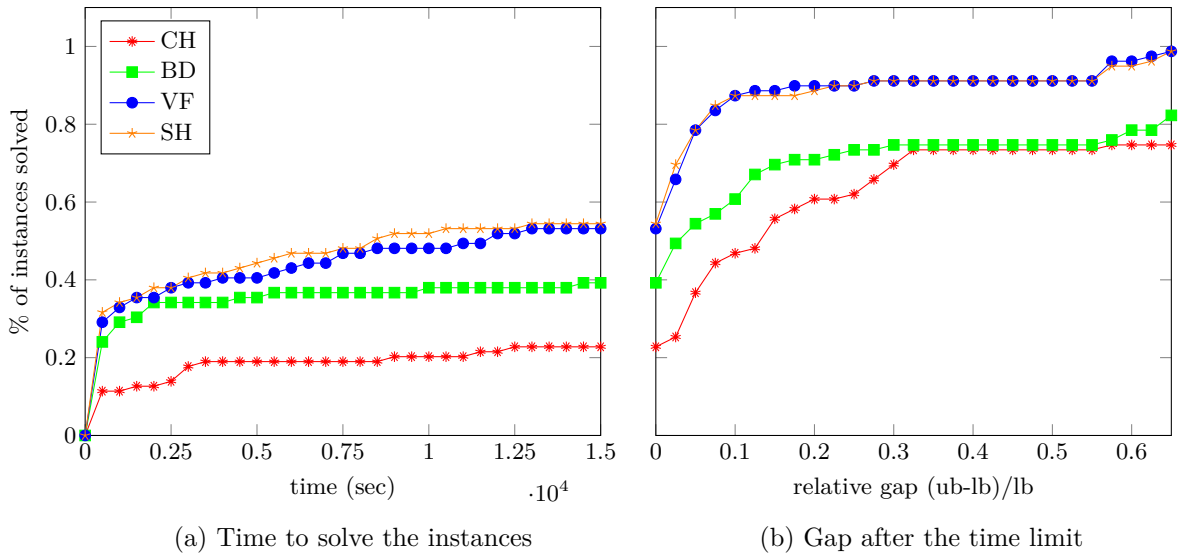

In this section we finally compare the B\&B algorithms using the
different upper bounds presented so far.
In Table~\ref{tab:nrBaBNodes} and Figure~\ref{fig:performance} we compare the
number of nodes generated in the B\&B algorithm as well as
the CPU time and the final gap.
The abbreviations refer to the following bounds used.

\begin{description}
  \item[\mytag{CH}] We consider the upper bound obtained by the ESCs in the 
  convex 
hull formulation~\mytag{CH} described in 
Section~\ref{sec:IEConvexCombination},
  \item[\mytag{BD}] solving this formulation 
with the bundle method \mytag{BD} as presented in 
Section~\ref{sec:IEConvexCombination}, 
  \item[\mytag{VF}] relaxing this formulation by considering 
only violated facets~\mytag{VF} as described in Section~\ref{sec:IEFacets}, 
  \item[\mytag{SH}] and using only separating hyperplanes~\mytag{SH} as 
  presented in 
Section~\ref{sec:IESeparatingHyperplanes}.  
  \item[\mytag{TH}] For better comparability we also consider our B\&B algorithm with only the 
Lovász theta function~\eqref{lovaszTheta} as an upper bound and denote this 
version with \mytag{TH}. Note that this boils down to 
solving~\eqref{SDPRelaxationWithESC} with $J=\emptyset$. 
\end{description}

If we are not able to solve an instance within the timelimit, we
indicate this in Table~\ref{tab:nrBaBNodes} by a cell that is colored 
{\textcolor{yellow!45}{\rule{0.8cm}{0.7\baselineskip}}}.
Whenever a cell is colored 
{\textcolor{red!25}{\rule{0.8cm}{0.7\baselineskip}}}
it means that the run did not finish 
correctly, for example because MOSEK produced an error or ran out of memory.  
Before discussing the results, we give the details on the instances as
well as on the soft- and hardware. 

\subsection{Benchmark Set and Experimental Setup}

We consider several different instances. First, we consider the
instances used 
in~\cite{GaarRendlFull}, i.e.,\ torus graphs, random near-$r$-regular 
graphs and random graphs from the Erd\H{o}s-Rényi model and also several 
instances from the literature.
Additionally we consider all instances from the DIMACS challenge for which the 
gap between $\vartheta(G)$ and $\alpha(G)$ is larger than one (i.e.,\ all 
instances which are not solved in the root node of our B\&B algorithm) and that 
have at most 500 vertices. 
Moreover, we consider spin graphs, which are produced with the command
\texttt{./rudy -spinglass3pm x x x 50 xxx1} 
(for $\texttt{x} \in \{5,7,9,11\}$)
by the graph generator 
rudy, which was written by Giovanni Rinaldi.\footnote{Available at 
    \url{www-user.tu-chemnitz.de/~helmberg/rudy.tar.gz}}

We implemented our B\&B algorithm with different upper bounds in \CPP. 
All computations were performed on an
Intel(R) Core(TM) i7-7700 CPU @ 3.60GHz with 32 GB RAM.
All programs were compiled with gcc version 5.4.0 with the optimization level 
-O3 and the CPU time was measured with \texttt{std::clock}. 
We set the 
random seed to zero. 
\elli{We use MOSEK~\cite{mosek} 8.1 in the methods \mytag{CH},
  \mytag{VF}, and \mytag{SH} to solve the
  SDPs~\eqref{SDPRelaxationWithESCConvHullOrig},
  \eqref{SDPRelaxationWithESCFacetInequalities}, and
  \eqref{SDPRelaxationWithESCSeparatingHyperplane}, respectively.
  Furthermore, we use it within the method \mytag{BD} for solving the
  subproblems within the bundle method for computing an
  approximate solution of~\eqref{SDPRelaxationWithESCConvHullOrig}. 
  Moreover, we use it to 
  solve the QP~\eqref{problemProjectionDistance} to compute the projection 
  distance
  in \mytag{SH} and when updating $J$, i.e., while searching for subgraphs with 
  violated ESCs and adding these subgraphs to $J$.
}
The execution time of our B\&B algorithm is limited to $4$ hours, i.e.,\ after 
this computation time we allow the B\&B to finish solving the SDP of the 
already started node in the B\&B tree and then stop.

{
\footnotesize
\centering
\begin{longtable}[htb]{|l|rrc|rrr|rr|}
\captionsetup{justification=centering}
\caption{The number of nodes in our B\&B 
algorithm}\label{tab:nrBaBNodes}\\
\hline 
Graphs&$n$ &$m$& $\alpha$ &\myref{TH} &\myref{CH} &\myref{BD} &\myref{VF} &\myref{SH} \\ 
\hline 
\endfirsthead
\caption{The number of nodes in our B\&B algorithm (cont.)}\\
\hline
Graphs&$n$ &$m$ & $\alpha$ &\myref{TH} &\myref{CH} &\myref{BD} &\myref{VF} &\myref{SH} \\ 
\hline 
\endhead
\hline \endfoot
ER\_100\_15& 100&   737& 24&     25   &      19   &      19   
&      19   &      19   \\ 
ER\_100\_25& 100&  1262& 17 &    27   &      49   &      23   &      21   
&      23   \\ 
ER\_100\_50& 100&  2501&  9 &   23   &      21   &      21   &      21   
&      21   \\ 
ER\_125\_15& 125&  1123& 27   & 243   &      22\nO&     173   &     139   &     
149   \\ 
ER\_125\_25& 125&  1928& 20   &  29   &      21   &      25   &      23   
&      21   \\ 
ER\_125\_50& 125&  3904& 9 &    69   &      63   &      63   &      63   
&      63   \\ 
ER\_150\_15& 150&  1664& 28 &   1469   &      41\nO&       -\nF&     835   &     
871   \\ 
ER\_150\_25& 150&  2810& 19  &   543   &      47\nO&     417   &     401   &     
429   \\ 
ER\_150\_50& 150&  5608& 10  &   79   &      58\nO&      77   &      77   
&      77   \\ 
ER\_175\_15& 175&  2213& 30 &  4633   &      35\nO&       -\nF&    1856\nO&    
2488\nO\\ 
ER\_175\_25& 175&  3777&  20 &  1191   &      32\nO&     818\nO&     985   &     
995   \\ 
ER\_175\_50& 175&  7578&11  &    91   &      24\nO&      87   &      87   
&      87   \\ 
ER\_200\_15& 200&  2917& 33 &   3149   &      21\nO&     479\nO&    1342\nO&    
1771\nO\\ 
ER\_200\_25& 200&  4924& 21 &   2883   &      18\nO&     274\nO&    1196\nO&    
1861\nO\\ 
ER\_200\_50& 200&  9913&  11 &  135   &      12\nO&      52\nO&     131   &     
131   \\ \hline
torus5&  25&    50& 10  &   3   &       1   &       1   &       1   &       1   
\\ 
torus7&  49&    98&  21  &  11   &       1   &       1   &       1   &       1   
\\ 
torus9&  81&   162&  36 &   33   &       1   &       7   &       1   &       1   
\\ 
torus11& 121&   242& 55  &   79   &       1\nO&      19   &       1   &       3   
\\ 
torus13& 169&   338& 78  &   211   &       6\nO&      59   &       5   &       5   
\\ 
torus15& 225&   450& 105  &   551   &       3\nO&       -\nF&      11   &      21   
\\ 
torus17& 289&   578& 136 &  1453   &       2\nO&       -\nF&      23   &      47   
\\ 
torus19& 361&   722& 171  &    -\nF&       -\nF&       -\nF&      31\nO&      
98\nO\\ 
torus21& 441&   882& 210  &  2956\nO&       -\nF&       -\nF&      23\nO&     
105\nO\\ 
torus23& 529&  1058& 253  & 1750\nO&       2\nO&      59\nO&      16\nO&      
89\nO\\ 
torus25& 625&  1250& 300  & 1140\nO&       -\nF&      47\nO&      10\nO&      
61\nO\\ 
torus27& 729&  1458& 351  &  774\nO&       -\nF&      33\nO&       7\nO&      
41\nO\\ 
torus29& 841&  1682& 406  &  520\nO&       -\nF&      23\nO&       5\nO&      
28\nO\\ 
torus31& 961&  1922& 465  &  369\nO&      12\nO&      17\nO&       4\nO&      
20\nO\\ 
torus33&1089&  2178& 528  &  261\nO&       -\nF&      12\nO&       2\nO&      
14\nO\\ 
torus35&1225&  2450& 595   & 188\nO&       -\nF&       9\nO&       2\nO&      
10\nO\\ 
torus37&1369&  2738& 666  &  134\nO&       -\nF&       6\nO&       2\nO&       
7\nO\\ 
\hline
reg\_n100\_r4& 100&   195& 40  &  145   &       8\nO&      71   &      11   &      
15   \\ 
reg\_n100\_r6& 100&   294& 34  &  251   &      32\nO&     123   &      47   &      
51   \\ 
reg\_n100\_r8& 100&   377& 31 &   141   &      18\nO&      85   &      37   &      
39   \\ 
reg\_n100\_r10& 100&   474& 28  &  147   &      45\nO&      95   &      51   
&      55   \\ 
reg\_n200\_r4& 200&   400& 81 &  9873   &      24\nO&       -\nF&     417\nO&     
647\nO\\ 
reg\_n200\_r6& 200&   593&$69 \leqslant \alpha \leqslant 72$ & 28221\nO&      24\nO&     642\nO&    1023\nO&    
1447\nO\\ 
reg\_n200\_r8& 200&   792&$60 \leqslant \alpha \leqslant 63$ &  27441\nO&      24\nO&     630\nO&    1242\nO&    
1663\nO\\ 
reg\_n200\_r10& 200&   980& $57 \leqslant \alpha \leqslant 59$ & 26591\nO&      24\nO&     606\nO&    1257\nO&    
1657\nO\\ \hline
rand\_n100\_p004& 100&   212&  45 &    3   &       1   &       1   &       1   
&       1   \\ 
rand\_n100\_p006& 100&   303& 38 &    23   &      11   &      13   &       5   
&       9   \\ 
rand\_n100\_p008& 100&   443& 32 &    35   &      11   &      19   &      11   
&      11   \\ 
rand\_n100\_p010& 100&   489&32  &    15   &       2\nO&       9   &       3   
&       3   \\ 
rand\_n200\_p003& 200&   631& 81  &   39   &       6\nO&      29   &       9   
&      11   \\ 
rand\_n200\_p004& 200&   816& 67 &   6165   &      23\nO&       -\nF&     
423\nO&     999\nO\\ 
rand\_n200\_p005& 200&   991& 64 &   879   &      19\nO&       -\nF&     279   
&     295   \\ \hline
CubicVT26\_5&  26&    39& 10  &    7   &       1   &       3   &       1   &       
1   \\ 
HoG\_6575&  45&   225& 10  &  189   &      79   &     105   &      77   &      
77   \\ 
Circulant47\_030&  47&   282& 13 &     5   &       1   &       1   &       1   
&       1   \\ 
PaleyGraph61&  61&   915& 5 &    49   &      47   &      47   &      47   &      
47   \\ \hline
spin5& 125&   375& 50 &   325   &       1   &       1   &       1   &       1   \\ 
spin7& 343&  1029&  $147 \leqslant \alpha \leqslant 151$ & 5992\nO&       -\nF&       -\nF&      79\nO&     181\nO\\ 
spin9& 729&  2187& $ 324 \leqslant \alpha \leqslant 333$  &  670\nO&       -\nF&      23\nO&       4\nO&      41\nO\\ 
spin11&1331&  3993& $ 605 \leqslant \alpha \leqslant 625 $  &  126\nO&       -\nF&       4\nO&       1\nO&       
8\nO\\ \hline
MANN\_a9&  45&    72&  16  &   7   &       3   &       5   &       3   &       3   
\\ 
MANN\_a27& 378&   702&  126  &   -\nF&       4\nO&     188\nO&     131\nO&     
580\nO\\ 
hamming6\_4&  64&  1312& 4 &    17   &       1   &       1   &       1   &       
1   \\ 
MANN\_a45&1035&  1980& 345 &    280\nO&       -\nF&      14\nO&      17\nO&      
43\nO\\ 
sanr200\_0\_9& 200&  2037& 42 &  5623   &      22\nO&     475\nO&    1102\nO&    
1693\nO\\ 
brock200\_1& 200&  5066& 21 &  1691   &      18\nO&     296\nO&    1073\nO&    
1379   \\ 
keller4& 171&  5100& 11  &   167   &      24\nO&     103\nO&     119   &     129   
\\ 
sanr200\_0\_7& 200&  6032& 18 &  1099   &      18\nO&     254\nO&     905\nO&     
935   \\ 
brock200\_4& 200&  6811& 17 &   709   &      16\nO&     185\nO&     619   &     
627   \\ 
brock200\_3& 200&  7852& 15   &  469   &      14\nO&     122\nO&     403   &     
409   \\ 
brock200\_2& 200& 10024& 12 &   127   &      14\nO&      54\nO&     123   &     
123   \\ 
c\_fat200\_5& 200& 11427& 58 &    11   &       1\nO&       7   &       7   &       
7   \\ 
p\_hat300\_3& 300&  11460& 36 &   622 \nO   &      4 \nO&       -\nF   &      58 \nO   &       
102 \nO  \\ 
brock400\_2& 400& 20014& 29  &   34\nO&       -\nF&       4\nO&      10\nO&      
14\nO\\ 
brock400\_4& 400& 20035& 33  &   34\nO&       -\nF&       4\nO&      10\nO&      
13\nO\\ 
brock400\_1& 400& 20077& 27  &   36\nO&       -\nF&       2\nO&       9\nO&      
12\nO\\ 
brock400\_3& 400& 20119& 31  &   35\nO&       -\nF&       2\nO&      10\nO&      
12\nO\\ 
p\_hat300\_2& 300&  22922& 25  &   24\nO   &      1 \nO&       -\nF  &       8\nO   &       
10 \nO  \\ 
sanr400\_0\_7& 400& 23931& 21 &    22\nO&       -\nF&       2\nO&       
6\nO&       8\nO\\ 
san400\_0\_7\_3& 400& 23940& 22  &    65\nO&       -\nF&       3\nO&       1   
&      14\nO\\ 
p\_hat500\_3& 500& 30950&  49  &   8\nO&       -\nF&       1\nO&       2\nO&       
4\nO\\ 
p\_hat300\_1& 300&  33917& 8 &   6 \nO   &      1 \nO&       -\nF  &       2\nO&       
2 \nO   \\ 
sanr400\_0\_5& 400& 39816& 13  &    5\nO&       4\nO&       1\nO&       
2\nO&       2\nO\\ 
\end{longtable}%
}

\subsection{First Computational Experiments}
\label{sec:firstComputationalComparison}

We first want to compare the two different versions of the B\&B
algorithm that use \myref{CH} and \myref{BD} 
to compute upper bounds,
i.e., we compare those versions that have already been established as 
upper bounds in the literature, but are now for the first time used within a 
B\&B algorithm.

First, by looking at Table~\ref{tab:nrBaBNodes} we observe that 20
instances were not solved correctly with \myref{CH},  
which is due to the fact that the SDPs to solve are huge and therefore MOSEK runs 
out of memory very often. Indeed, by using \myref{BD} and hence not having to 
solve that large SDPs only 10 instances are not solved correctly, most of 
them due to other MOSEK errors.
When we compare the number of B\&B nodes for the instances where both 
\myref{CH} and \myref{BD} finished we see that typically the number is the 
same, or there are slightly more nodes for \myref{BD}. This is plausible, as we 
only have an approximate solution of~\eqref{SDPRelaxationWithESC} when 
using \myref{BD}, but an exact solution of~\eqref{SDPRelaxationWithESC}  with 
\myref{CH}.

We next take a closer look on the lines labeled \myref{CH} and
\myref{BD} in the performance profiles in Figure~\ref{fig:performance}.
The B\&B code using \myref{BD} as bounding routine can solve much more instances
within a given time than when \myref{CH} is used. 
Once the time limit is reached, the gap is typically much lower for
\myref{BD} than it is for \myref{CH}.

In a nutshell, even though \myref{BD} solves only a relaxation 
of~\eqref{SDPRelaxationWithESC}, using it is much faster than using \myref{CH} while it 
does not increase the number of B\&B nodes a lot. This justifies 
considering only a relaxation of~\eqref{SDPRelaxationWithESC}. 

\subsection{Computational Experiments with New Relaxations}

Up to now we have used the exact subgraph approach of~\cite{AARW} with the 
implementation proposed by Gaar and Rendl~\cite{GaarRendlIPCO,GaarRendlFull} in 
order 
to get tight upper bounds on the stability number within a B\&B 
algorithm for 
solving the stable set problem. 
So far we have proven the strength of the bounds by showing that the number of 
nodes in a B\&B algorithm reduces drastically by using these 
bounds, 
however the computational costs are enormous in the original version with 
the convex hull formulation  
\myref{CH} and they are still substantial with the bundle approach 
\myref{BD} from~\cite{GaarRendlIPCO,GaarRendlFull}.
Therefore, we now discuss the numerical results of the B\&B algorithm
using the new relaxations \myref{VF} and \myref{SH}.

%

Looking at Table~\ref{tab:nrBaBNodes},
the first thing we observe is that both~\myref{VF} and~\myref{SH} never lead to 
a MOSEK error, hence they are more robust than the other versions, presumably 
due to their smaller size of the SDPs to solve in the B\&B nodes.
For 9 instances both~\myref{VF} and~\myref{SH} and in 13 instances at least one 
of~\myref{VF} and~\myref{SH} is able to finish within the time limit for an 
instance where both~\myref{CH} and~\myref{BD} were not able to finish. 

If we compare the number of B\&B nodes in Table~\ref{tab:nrBaBNodes} for the 
finished instances, then we see 
that typically the number of B\&B nodes for~\myref{VF} is smaller than 
those of~\myref{BD}, which makes sense because in~\myref{BD} we only 
approximately solve~\eqref{SDPRelaxationWithESC} whereas in~\myref{VF} we solve 
a possibly very tight relaxation of~\eqref{SDPRelaxationWithESC} exactly. The 
number of B\&B nodes needed by~\myref{SH} is typically a little bit larger than 
the one of~\myref{VF}, which is in tune with the 
\elli{
empirical finding that~\eqref{SDPRelaxationWithESCFacetInequalities} gives 
better bounds than~\eqref{SDPRelaxationWithESCSeparatingHyperplane} in the 
small example considered in Section~\ref{sec:SmallComputationalComparisonVFSH}.}
In a nutshell, for finished runs typically~\myref{CH} and~\myref{VF} need 
roughly the same number of nodes,~\myref{SH} needs a little bit and~\myref{BD} 
needs many more nodes in the B\&B tree.

As for the running times, in Figure~\ref{fig:performance} we see that both \myref{VF} and \myref{SH} are faster than
\myref{BD} and considerably faster than \myref{CH}.
\myref{VF} is a little bit slower than~\myref{SH}.
For those instances that cannot be solved to optimality, the gap when
the time limit is reached is roughly the same for \myref{VF} and
\myref{SH}, and it is considerably tighter than for \myref{BD}.


We have demonstrated that within a B\&B algorithm both our 
relaxations~\myref{VF} and~\myref{SH} work better than already existing 
SDP based methods. 
In particular using~\myref{SH} allows to keep the majority of the 
strength of the upper bound~\eqref{SDPRelaxationWithESC} (i.e.,\ keeping the 
number of vertices in the B\&B tree low) by reducing the running time
so that within the time limit almost 60\% of the instances are solved,
as compared to \myref{CH} that only manages to solve a bit more than
20\%.

\section{Conclusions}
\label{section:conclusions}

We introduced an algorithm for computing the stability number of a
graph using semidefinite programming.
While there exist several exact solution methods for finding the
stability number, those based on semidefinite programming are rather rare.

We implemented a B\&B algorithm using the SDP relaxations introduced
in~\cite{GaarRendlIPCO,GaarRendlFull} together with heuristics from the literature. Moreover, we
further relaxed the SDPs, getting more tractable SDPs still producing
high-quality upper bounds. This is confirmed by the numerical
experiments where we compare the number of nodes to be explored in the
B\&B tree as well as the CPU times.

While SDPs produce strong bounds, the computational expense for
solving the SDPs is sometimes huge. 
In particular, there is potential for improvement of the running time for 
solving SDPs with many constraints.
We use MOSEK as a solver, which uses the interior point method to solve an SDP.
For large instances it would be beneficial to use a solver based on the boundary 
point method~\cite{bpm06,bpm09} or DADAL~\cite{dadalPaper}.
Moreover, the solver computing $\vartheta^+$ as an upper bound~\cite{cerulli2020}
combined with the relaxations above, \elli{may} push the performance of the
B\&B solver even further.
Also, these other solvers are capable of doing warm starts, that can
have big advantages within a B\&B framework.
Since all these implementations are available as MATLAB source code
only, they need 
to be translated to C or \CPP\  first. This will be part of our future
work.

Another line of future research is working out different strategies
for identifying violated subgraphs, that should also lead to a more
efficient overall algorithm. 



\ifSpringer
\else
\bibliographystyle{plain}
\fi
\bibliography{diplomarbeit}

\ifSpringer
\else
{\small
  \vspace*{1ex}\noindent
  Elisabeth Gaar,
  \href{mailto:elisabeth.gaar@jku.at}{\url{elisabeth.gaar@jku.at}},
  Johannes Kepler University Linz,
  Altenberger Stra{\ss}e~69, 4040 Linz, Austria

  \vspace*{1ex}\noindent
  Melanie Siebenhofer,
  \href{mailto:melanie.siebenhofer@aau.at}{\url{melanie.siebenhofer@aau.at}},
  Alpen-Adria-Universität Klagenfurt,
  Universitätsstraße 65--67, 9020 Klagenfurt, Austria
  
  \vspace*{1ex}\noindent
  Angelika Wiegele,
  \href{mailto:angelika.wiegele@aau.at}{\url{angelika.wiegele@aau.at}},
  Alpen-Adria-Universität Klagenfurt,
  Universitätsstraße 65--67, 9020 Klagenfurt, Austria
}
\fi

\end{document}